\input amstex

\def\bb#1#2{b^{#1}_{#2}}

\def\bcc#1{\Bbb C^{#1}}
\def\bll#1{\Bbb L^{#1}}
\def\brr#1{\Bbb R^{#1}}
\def\brrr{\Bbb R }

\def\R{\Bbb R}
\def\bpp#1{\Bbb P^{#1}}
\def\bccc{\Bbb C }

\def\cf{\Cal F}

\def\ee#1{e_{#1}}
\def\uee#1{e^{#1}}

\def\ue#1{e^{#1}}

\def\gg#1#2{g_{{#1}{#2}}}
\def\ugg#1#2{g^{{#1}{#2}}}

\def\hh#1#2#3{h^{#1}_{{#2}{#3}}}
\def\hd{, \hdots ,}
\def\inv{{}^{-1}}

\def\k{\kappa}

\def\lm#1{\lambda_{#1}}
\def\ooo#1#2{\omega^{#1}_{#2}}
\def\oo#1{\omega^{#1}}
\def\ot{\!\otimes\!}
\def\pp#1#2{\pi_{{#1}{#2}}}

\def\qq#1#2#3{q^{#1}_{{#2} {#3}}}

\def\ra{\rightarrow}
\def\rr#1#2#3#4{r^{#1}_{{#2} {#3}{#4}}}

\def\Sum#1{\underset{#1}\to\Sigma}

\def\tee#1{\tilde e_{#1}}
\def\tg#1#2{\tilde g_{{#1}{#2}}}
\def\tmu#1#2{\tilde \mu^{#1}_{#2}}
\def\tm{M} 
\def\tnu#1#2{\tilde \nu^{#1}_{#2}}

\def\tcf{\tilde\Cal F}
\def\tpi{\tilde\pi}

\def\tooo#1#2{\tilde\omega^{#1}_{#2}}
\def\to#1{\tilde\omega^{#1}}
\def\too#1#2{\tilde\omega^{#1}_{#2}}
\def\tO#1#2{\tilde\Omega^{#1}_{#2}}

\def\tdim{\text{dim}}

\def\tmod{\text{ mod }}

\def\upperp{{}^{\perp}}

\def\ww{\wedge}
\def\xx#1{x^{#1}}

\define\intprod{\mathbin{\hbox{\vrule height .5pt width 3.5pt depth 0pt %
        \vrule height 6pt width .5pt depth 0pt}}}

\newcount\nummereq
\nummereq=0
\def\labeleq#1{\global\advance\nummereq
by1\global\edef\here{\number\nummereq}\here
\global\let #1 =\here}

\def\strutdepth{\dp\strutbox}
\def\leftnote#1{\vadjust{\vbox to 0pt{
   \vss\hbox to\hsize{\hskip-0.9in
   \vbox{\hsize0.8truein\overfullrule0pt\noindent
       {\it #1}}\hss}\vskip\strutdepth}}}
\def\linebreak{\hfill\break}

\documentstyle{amsppt}
\magnification = 1200
\hsize =15truecm
\hcorrection{.5truein}
\baselineskip =18truept
\vsize =22truecm
\NoBlackBoxes
\topmatter
\title
On isometric and minimal isometric embeddings \endtitle
\author
Thomas Ivey and J.M. Landsberg
\endauthor
\thanks  The second author is partially supported by  NSF grant DMS-9303704.
\endthanks
\address {Dept. of Mathematics, Case Western Reserve Univ.,
Cleveland OH 44106-7058}\endaddress
\email {txi4\@po.cwru.edu} \endemail
\address{Department of Mathematics,
Columbia University, New York,  NY 10027}
\endaddress
\curraddr Math\'ematique,
Laboratoire Emile Picard,
Universit\'e Paul Sabatier,
31062 Toulouse CEDEX,
France  \endcurraddr
\email {jml\@picard.ups-tlse.fr }
\endemail
\date {  July 31, 1996 }\enddate
\endkeywords
\abstract{ In this paper we study  critial isometric and minimal isometric
embeddings of classes of Riemannian metrics
 which we call {\it quasi-$\k$-curved metrics}. Quasi-$\k$-curved
metrics generalize the metrics of space forms. We construct
  explicit examples and prove results about  existence and rigidity.
}\endabstract

\endtopmatter

\document
\heading Introduction \endheading

\noindent{\bf Definition}: Let $(M^n,\tilde g)$ be a Riemannian manifold.
We will say $\tilde g$ is a
{\it quasi-$\k$-curved metric} if there exists a
smooth positive definite quadratic form $Q$ on $M$
such that for all $x\in M$
$$
R_x=
-\gamma (Q_x,Q_x) +(\k+1)\gamma(\tilde g_x,\tilde g_x) \tag\labeleq{\quasik}
$$
where $\gamma :S^2T^*\ra S^2(\Lambda^2T^*)$
denotes the algebraic Gauss mapping and
$R_x$ the Riemann curvature tensor. (See \S 1. for more details.)

\smallpagebreak

Quasi-$\k$-curved  metrics are a generalization of
{\it quasi-hyperbolic metrics} defined in [BBG], which
correspond to $\k=-1$. We will also refer to the
case $\k=0$ as {\it quasi-flat metrics}. We will assume that $n\ge 3$. When
$n = 3$, the quasi-$\k$-curved condition is an open condition on the metric,
and thus in this case the class of metrics we study is quite general.
The condition  is stronger in higher dimensions.

In this paper we study local isometric embeddings and
minimal isometric embeddings of     quasi-$\k$-curved manifolds.
Before giving our results, it will be useful to review some of what is known:

\smallpagebreak

\noindent{\bf Local isometric embeddings}

Given a Riemannian manifold $(M^n,\tilde g)$, one may ask if
it admits a local isometric embedding into a Euclidean space
$\brr{n+r}$ or more generally a space form $X(\epsilon)^{n+r}$ of constant
sectional curvature $\epsilon$.
If $M$ is more positively curved than $X$, one expects to
have local isometric embeddings (e.g. the embedding of $S^n$ into
$\brr{n+1}$).
We will be concerned with the case $M$ is less positively curved than $X$.

The critical codimension for the
isometric embedding problem is $r=\binom n2$. The
Cartan-Janet theorem states that local isometric embeddings of
analytic Riemannian $n$-folds into $\brr{n+\binom n2}$ exist and depend locally
on a choice of $n$ arbitrary functions of $n-1$ variables.
(These ``dimension counts'' come from the Cartan-K\"ahler Theorem [Car2].)
In this paper, we will be interested in the overdetermined case $r<\binom n2$.

In the most overdetermined case ($r=1$), Thomas [T] observed that
the Codazzi equations of a hypersurface with non-degenerate
second fundamental form are consequences of the Gauss equations
when $\tdim (M)\geq 4$. Thus codimension one questions reduce
to questions in multi-linear algebra.
(See [CO] for a clear exposition.)

For $r \le n$, Cartan [Car1] studied the isometric immersions of a flat space
into
a Euclidean space, showing that if $r\leq n-1$ there is no local
isometric embedding (other than the totally geodesic one) and
that when $r=n$ such embeddings depend on $\binom n2$ functions of $2$
variables.
In the course of his proof, Cartan proved a basic theorem about exteriorly
orthogonal symmetric bilinear forms (see e.g. [Spivak], V.11.5).
Using a reduction  to the flat case (see \S 1), Cartan went on to show that
hyperbolic space $H^n$ admits no local isometric embedding into $\brr {2n-2}$
but admits local isometric embeddings into $\brr{2n-1}$.  Moreover, these
embeddings depend on a choice of $n^2-n$ functions of one variable,
which is the most possible for local
isometric embeddings (with nondegenerate second fundamental form)
for any $n$-fold into  $\brr{2n-1}$.

  Cartan's work implies that for such
isometric embeddings of hyperbolic space there is, at each point,
an orthonormal basis of the tangent space under which
the second fundamental form $II$ is diagonalized with
respect to the metric. Thus, for such embeddings there are principal
tangent directions and  princial normal vectors, in analogy with the case
of hypersurfaces.

Quasi-$\k$-curved metrics share the property that,
for immersions to $X(\k+1)$ in the critical
codimension, the tangent space admits a
basis for which the second fundamental form  is diagonalized,
but the basis is $Q$-orthonormal instead
of $\tilde g$ orthonormal.  We will   refer to these vectors
 as the {\it principal tangents}, and the corresponding values of $II$ as
 {\it principal normal vectors}.

Cartan's results for hyperbolic space were generalized by
Berger, Bryant and Griffiths [BBG]
to quasi-hyperbolic metrics. They showed that no quasi-hyperbolic metric
 admits a local isometric embedding
to Euclidean space when $r<n-1$, and characterized the quasi-hyperbolic
metrics that admit local isometric embeddings in the
critical case $r=n-1$.

To understand this characterization,
note that the Riemann curvature tensor for a quasi-hyperbolic
metric is like that of a hypersurface in Euclidean space but the
sign is wrong.  If instead one considers
spacelike hypersurfaces in Lorentz space, one obtains
quasi-hyperbolic metrics where $Q=II$. In [BBG] they assert that {\it all}
 quasi-hyperbolic metrics arise in this way (at least when $n>3$).

The sub-class of
quasi-hyperbolic metrics  satisfying the integrability conditions
to admit a critical isometric embedding are precisely the spacelike
hypersurfaces satisfying the condition
$$\nabla II_{\Bbb L}= L\cdot  II_{\Bbb L}, \tag\labeleq{\nabtwo} $$
where $II_{\Bbb L}$ denotes the second fundamental form
of the hypersurface embedding,
$\nabla$ is the connection form of
the metric $\tilde g$ and $L \in \Omega^1(M)$ is some linear form.
(This is in fact an intrinsic condition because $II_{\Bbb L}=Q$.)
A result in projective geometry ([GH], B.16),
which specializes to affine geometry, shows
that (\nabtwo) occurs if and only if $\tm^n$ is embedded as a patch
of a quadric hypersurface in affine space.
In [BBG] they also assert that in the $n=3$ case all quasi-hyperbolic
$3$-folds satisfying the integrability conditions occur as quadric
hypersurfaces in $\bll 4$.
Metrics satisfying (\nabtwo) admit
local isometric embeddings into $\brr{2n-1}$ which
depend on a choice of $n^2-n$ functions of one variable,
as in the case of hyperbolic space.

The [BBG] result generalizes to quasi-$\k$-curved metrics:
\proclaim{Theorem A}  Let $(M^n,g)$, $n\ge 3$, be a quasi-$\k$-curved
 Riemannian manifold.
Let $X^{2n-1}(\k+1)$ be a space form with constant sectional
curvature $\k+1$. Then there exist local isometric embeddings
$M^n \hookrightarrow X^{2n-1}(\k+1)$, with local solutions depending
on $n^2-n$ functions of one variable, if and only if $\nabla Q$
is a symmetric cubic form on $M$ and
$$\nabla Q = L\cdot Q$$
for some linear form $L\in \Omega^1(M)$.
\endproclaim

To produce quasi-hyperbolic metrics, it was natural to
look at hypersurfaces in $\bll{n+1}$. To study
quasi-$\k$-curved metrics, it is natural to look for codimension
two spacelike submanifolds of $\bll{n+2}$ having a spacelike
section $\sigma$ of the unit normal bundle such that
$\sigma\intprod II = \sqrt{\k+1}g$. Such $M$ are quasi-$\k$-curved
with quadratic form $Q= \sigma\upperp\intprod II$, where $\sigma\upperp$
is the timelike section of length $-1$ normal to $\sigma$.  (We do not
know whether or not all quasi-$\k$-curved metrics arise in this way.)
Using methods developed in [Lan], we first show any such
spacelike submanifold
is in fact a hypersurface inside a sphere:

\proclaim{Lemma} Let $M^n\subset \bll{n+2}$ be a spacelike submanifold
of Lorentz space.
 Suppose there exists a section $\sigma$ of the unit normal bundle of $M$,
such that
$$
\sigma\intprod II = cg,
$$
where $c$ is some constant and $g$ is the induced metric.
Then   $M$ is conguent to
a submanifold of the Lorentzian sphere of radius $1/c$.
\endproclaim

\noindent{\bf Remark}: The lemma is valid in more general contexts
which will be suggested in the proof.

\proclaim{Corollary}Spacelike hypersurfaces of the Lorentzian sphere of
radius $\frac 1{\sqrt{\k+1}}$
in $\bll{n+2}$ are quasi-$\k$-curved.
\endproclaim

Let $x^0, x^1, \cdots, x^{n+1}$ be coordinates on $\bll{n+2}$
such that
$$<x,x> = -(x^0)^2 + (x^1)^2 + \cdots + (x^{n+1})^2.$$
Then the sphere our computations will produce has the equation:
$$
-(x^0)^2 + (x^1)^2 + \cdots + (\xx n)^2 +(x^{n+1}-\frac 1{2\sqrt{\k+1}})^2
=\frac 1{ \k+1}\tag\labeleq{\quadric}
$$
with $\sigma_{(0\hd 0)}=\partial/\partial\xx{n+1}$ and
$T^*_{(0\hd 0)}M=\{ d\xx i \}$, $ 1\leq i\leq n$.
Note that for the quasi-hyperbolic case,
(\quadric) specializes to the linear
subspace $\xx{n+1}=0$; the sphere of radius infinity.

We also classify those hypersurfaces that satisfy the isometric embedding
criteria, depending on the multiplicity of the eigenvalues of $Q$ with
respect to $g$.
\proclaim{Theorem B} Let $M^n\subset \bll{n+2}$ be quasi-$\k$-curved
with nondegenerate quadratic form $Q$.

Case 1. If, on an open set of $M$,
there exists an eigenvalue of $Q$ of multiplicity one, then
$\nabla Q = LQ$ for some linear form $L$ if and only if
$M$ is the intersection of (\quadric) and the quadric
$$
0=\xx 0 - q_{ij}\xx i\xx j - \lambda_j\xx 0\xx j - b(\xx 0 )^2
-(\k +1)\xx 0\xx{n+1}
\tag\labeleq{\tquadric}
$$
where $q_{ij}, \lambda_j$ are respectively the coefficients
of $Q$ and $L$
at $(0\hd 0)$ with respect to the orthonormal basis $d\xx i$,
and $b$ is a constant.

Case 2. If, on an open set of $M$,  there are no eigenvalues
of muliplicity one, then $\nabla Q=LQ$ if and only if
$\nabla Q=0$, in which case $M$ is a product of space forms and is the
intersection
of (\quadric) with  the quadric
$$
0=\xx 0 - q_{ij}\xx i\xx j- b(\xx 0 )^2-(\k +1)\xx 0\xx{n+1}
\tag\labeleq{\zquadric}
$$
where $q_{ij}$  are   the coefficients
of $Q$
at $(0\hd 0)$ with respect to the orthonormal basis $d\xx i$,
and $b$ is a constant.
\endproclaim

\bigpagebreak

\noindent{\bf Minimal isometric embeddings}

\smallpagebreak

Thanks to the work of Calabi,
there is
a reasonable understanding of the Riemannian metrics of
minimal surfaces. However
in higher dimensions, almost nothing is known
beyond   algebraic  restrictions
coming from the  Gauss equations
and restrictions coming from the isometric embedding alone.

Calabi showed that any
surface $M^2\subset \brr n$ that is minimally and
isometrically embedded must arise as a projection from some
Hermitian isometric embedding
of $M$ into some $\bcc m$, where  $m\leq n\leq 2m$ (see [Law]).  From this
description, it follows that
any such $M$ always admits some number of constants' worth (again, in the
Cartan-K\"ahler language) of
noncongruent
minimal isometric embeddings,
but never any functions' worth. For example,
non-planar minimal surfaces
in $\brr 3$ always admit a one parameter
family of minimal isometric deformations,
the most famous of which is the family connecting catenoids and helicoids.
One may wish to contrast this situation with the isometric embedding
problem for surfaces in $\brr 3$, without the requirement of minimality,
where, as stated above, a generic metric admits
two functions of one variable's worth of noncongruent isometric embeddings
([Car1]). (In fact, one can pose a Cauchy problem with a space curve as
initial data and realize these two functions as the
curvature and torsion of the curve.)

It follows from Calabi's work that no patch of the hyperbolic
plane admits a minimal isometric embedding to any
finite dimensional Euclidean space. On the other hand, Calabi
showed that all finite dimensional hyperbolic spaces admit
a minimal isometric embedding into a Hilbert space ([Cal]).

Moore [M] proved that the only ways in which an $n$-dimensional space
form $M^n$ can be locally isometrically embedded as a minimal submanifold in
a   space form $X^N(\epsilon )$  with $N\leq 2n-1$
are if the image is totally geodesic, or
$M^n$ is flat and its image is a piece of the Clifford torus $T^n$
 in $S^{2n-1}$.
In particular, $H^n$ admits no local minimal
isometric embeddings into $\R^{2n-1}$.

Define the {\it rank} of an embedding as
the dimension of the image of the second fundamental
form as a
linear map
$$\text{II} \, : \, S^2 T_x M \longrightarrow N_x M.$$
In other words, the dimension of
the second osculating space of $M$ at a point is $n$ plus the rank.
We obtain the following extensions of Moore's theorem:
\proclaim{Theorem C}

1. If $(M^n,\tilde g)$
 is a quasi-$\k$-curved
manifold and $Q$ is $\tilde g$-parallel, then
$M$ does not admit any local embedding of rank less
 than $n$ into $X^N(\epsilon)$, with $\epsilon < \k$,
which is isometric and minimal,
or which is isometric with parallel mean curvature vector,
except when $M$ is flat and its image is a piece of the Clifford torus.

2. If an $n$-dimensional space
form $M^n(\k )$ is locally isometrically embedded as a minimal submanifold
 of constant rank in an $N$-dimensional space form $X^N(\epsilon)$,
then either $M$ is totally geodesic, or the rank is at least $n-1$.  If the
rank is exactly $n-1$, then either $M$ is totally geodesic or $M$ is flat
and its image is a piece of the Clifford torus $T^n$ in $S^{2n-1}$.
(In particular, $H^n$ admits no local minimal
isometric embeddings of constant  rank $n-1$ into $\R^{N}$.)
\endproclaim

We also obtain the following rigidity theorem:
\proclaim{Theorem D} Let $M^n$ be a quasi-$\k$-curved Riemannian
manifold.  Let $\ell$ be the dimension of the principal orbits of the
identity component of the isotropy group of $M$, at a generic point
$p \in M$, acting on the $Q$-orthonormal frames of $T_pM$.
Then the minimal isometric embeddings,
and more generally the isometric embeddings with parallel
mean curvature vector, of $M$ into
a space form $X^{2n-1}(\k+1)$ depend, up to rigid motions, on at most a choice
of ${n \choose 2}-\ell$ constants.
\endproclaim
(Note that, when $\gamma$ is followed by the Ricci trace, the result is a
map of full rank on the space of quadratic forms at a positive definite $Q$.
Hence the identity component of the isometry group of $M$ preserves $Q$.)

The minimality condition imposes additional integrability conditions on the
metric; these take the form of an overdetermined
set of polynomials involving $Q$, the curvature
tensor, and their covariant derivatives.  Because of Theorem C, we know these
additional conditions are non-trivial, and so we obtain

\proclaim{Theorem E} A  non-empty Zariski-open subset of the  space of
quasi-$\k$-curved manifolds $M^n$
that admit a local isometric embedding into a space
form $X^{2n-1}(\k+1)$  does not admit any minimal isometric embedding.
\endproclaim
\smallpagebreak
\noindent{\bf Outline of the paper}
\medskip
In \S 1 we review the algebraic form of the Gauss equations of a submanifold
of a space form and explain the  quasi-$\k$-curved condition in more detail.
In \S 2 we set up the isometric embedding problem following [BBG]
and prove Theorem A. In
\S 3 we describe what happens to the system when we add the
 minimality condition, and we prove Theorems C and    D
in the case of minmality.  In \S 4 these results
are extended to the case of parallel mean curvature vector.
In \S 5 we construct
quasi-$\k$-curved $n$-folds as submanifolds of
 spheres in $\bll{n+2}$, proving the   Lemma and Theorem B.

\heading \S1. The Gauss equations \endheading
The class of Riemannian metrics we will be dealing with have a special
property that is best described in terms of the algebraic form of the
Gauss equations of a submanifold of a space form.

Throughout this paper, let $V=\brr n$ and $W=\brr r$; we endow $W$ with
the standard inner product.

Let $K $ be the kernel of
the skew-symmetrization map:
$$0 \longrightarrow K \longrightarrow
 \Lambda^2 V^*\ot\Lambda^2 V^*
\longrightarrow V^* \otimes \Lambda^3 V^*. $$
Note that actually $K\subset S^2(\Lambda^2 V^*)$. (This is
often called the first Bianchi identity.)
$K$ is the space of tensors with the symmetries of the Riemann curvature
tensor.

Let
$$\gamma:S^2V^*\ot W \ra K$$
be the $Gl(n)\times O(r)$-equivariant quadratic map defined as follows:
in terms of an arbitrary basis $\{\uee i \}$ for $V^*$ and
an orthonormal basis $\{ w_{\mu}\}$ for $W$,
$\gamma$ is given by
$$\hh {\mu} ij \ue i\ue j\ot w_{\mu}\mapsto
\Sigma_{\mu}(\hh {\mu} ik\hh {\mu} jl -\hh {\mu} il\hh {\mu} jk)
(\ue i\ww\ue j)\otimes (\ue k\ww\ue l).$$
For a submanifold $M^n \subset X^{n+r}(\epsilon )$ with the induced metric,
we will have $V=T_xM$, $W=N_xM$ (the fibre of the normal bundle at $x$),
and
$$R=\gamma (II,II)+ \epsilon\gamma (g,g)\tag\labeleq{\gauss}$$
where here and in what follows, to apply $\gamma$ to an
element of $S^2V^*$, just take $W=\brrr$.

\smallpagebreak

As mentioned
in the introduction,
Cartan realized that one could study isometric embeddings of hyperbolic
space via the isometric embeddings of flat space.  For, the curvature
tensor $R_0$
of hyperbolic space has the property
that it is minus the image of
the metric   $g_0 \in S^2 V^*$ under the map $\gamma$:
$$
R_0 = -\gamma(g_0,g_0).\tag\labeleq{\aquasi}
$$
Letting $\hat W = W\oplus \brrr$ and
$\widehat{II} = II \oplus g_0$, in the case of hyperbolic space, one obtains
$$0= \gamma (\widehat{II},\widehat{II} ),\tag\labeleq{\flatgauss}$$
the   Gauss equations
for an isometric embedding of a flat metric into $\brr{n+r+1}$.

The Gauss equations for embedding a quasi-$\k$-curved metric $\tilde g$
into a space form of curvature $\k+1$ have the form
$$-\gamma (Q,Q) +(\k+1)\gamma(\tilde g,\tilde g)
= \gamma (II,II)+ (\k+1)\gamma (\tilde g,\tilde g)$$
When one defines $\widehat{II} = II \oplus Q$, once again the Gauss equations
take the form (\flatgauss).

\heading \S2. Moving frames and the isometric embedding system\endheading
In this section we summarize what we will need from [BBG].

Let $(\tm,\tilde g)$ be a Riemannian manifold.
Let $\tpi  :\tcf \ra \tm$ denote the bundle of  all
frames of $T\tm$, i.e.,   the fiber over a point
$x\in \tm$ is the set of all bases $\tee 1\hd\tee n$  of $T_x\tm$.
On $\tcf$,
write $dx = \to i\tee i$ where
the one-forms $\to i$ are semi-basic to the projection
$\tpi$.
Let $\tg ij = \tilde g(\tee i, \tee j)$.  On $\tcf$ we have structure equations
$$
\align
& d\tg ij = \tg ik\tooo kj + \tg kj\tooo ki \\
&d\to i = -\tooo ij\ww\to j\\
&d\tooo ij + \tooo ik\ww\tooo kj = \tO ij,
\endalign
$$
where  the forms $\tooo ij$ are connection forms for the
Levi-Civita connection associated to $\tilde g$
and the forms $\tO ij$ are the curvature two-forms for this connection.

We will  set up an isometric embedding
system for quasi-$\k$-curved metrics following [BBG].
Let $\epsilon=\k+1$ and let $X(\epsilon)^{n+r}$ denote
the space form of constant sectional curvature $\epsilon$.
Let $\cf\ra X(\epsilon)^{n+r}$ denoted the frame bundle adapted
such that
for $x\in X$ the fibre of $\cf$ consists of bases
$(\ee i, \ee\mu)$ for $T_xX$, such that
$$
\ee i \cdot \ee \mu =0,\qquad\ee\mu\cdot \ee\nu = \delta_{\mu\nu}
$$
(We will use index ranges $1\leq i,j,k\leq n$ and $n+1\leq \mu,\nu\leq n+r$.
The ``$\cdot$'' is the standard inner product.)
Let $g_{ij}=\ee i\cdot \ee j$
and let $g^{ij}$ be the components of $g\inv$.  If we treat $f = (x,\ee i,
\ee\mu)$ as a matrix-valued function on $\cf$, and write $df=f \Omega$, then
$\Omega$ is a matrix-valued 1-form, the {\it Maurer-Cartan form}.
We will denote the entries of $\Omega$ as follows:
$$
\Omega =
\pmatrix
0&-\epsilon \oo j&-\epsilon\oo\nu\\
\oo i&\ooo ij & \ooo i\nu \\
 \oo\mu &\ooo \mu j & \ooo\mu\nu\endpmatrix
$$
where we have symmetries
$$\ooo\mu\nu = -\ooo\nu\mu,\qquad
\ooo k\nu =-\ugg ki\ooo\nu i.
$$
It follows from the definition of $\Omega$ that
$$
d\gg ij = \gg ik\ooo kj + \gg jk \ooo ki.
$$
The {\it Maurer-Cartan equation},
 $d\Omega = -\Omega\ww\Omega$
enables us to compute the exterior derivatives of the components of $\Omega$.

On the submanifold $\Sigma\subset \cf\times\tilde\cf\times (S^2V^*\ot W)$
defined by the equations $\tg ij - \gg ij =0$,
   define the Pfaffian system
$$
I^{std} =
\{ \oo i-\to i ,\oo\mu ,\ooo ij -\tooo ij ,\ooo\mu i-\hh\mu ij \oo j\}.
$$
with independence condition
$\to 1\ww\hdots\ww\to n \ne 0$.
Integral $n$-manifolds of this system
 are graphs, on the level of frames, of
isometric embeddings of $(\tm,\tilde g)$.
By differentiating the last set of forms in $I^{std}$,
one sees that integral manifolds can
only lie in the subset $\Sigma' \subset \Sigma$ where the Gauss equations
(\gauss) are satisfied. Also, note that this EDS is invariant under the group
$GL(n) \times O(r)$, acting by orthonormal changes of basis among
the $\ee \mu$ and
arbitrary but simultaneous changes of basis among $\tee 1\hd \tee n$ and
among $\ee 1\hd\ee n$.

Henceforth we will assume that $r=n-1$ and that the metric $\tilde g$
 is quasi-$\k$
curved with respect to a positive definite quadratic form $Q$.

Now the Gauss equations take the form (\flatgauss).  By an application of
Cartan's theorem on exteriorly orthogonal forms, there exists at each point
a $Q$-orthonormal basis for $TM$ such that:
$$
\hh\mu ij = \delta_{ij}\bb\mu i,\tag\labeleq{\hform}
$$
and the $n$ vectors $b_i = (b_i^\mu)\in W$ satisfy
$$b_i \cdot b_j = -1\ \text{for}\ i\ne j.\tag\labeleq{\bcond}
$$
(See [BBG], \S\S 4.2ff for details.) Note that such a basis is unique
up to permutation.

\smallpagebreak

Let $\pi_Q :\tcf_Q\ra \tm$ denote the bundle of $Q$-orthonormal
frames of $TM$, to which we restrict all the forms defined on $\tcf$.

Let $\Cal W$ be the smooth submanifold of $V^* \otimes W$ defined by (\bcond).
On the submanifold $\Sigma'' \subset \tcf_Q\times\cf\times \Cal W$
defined by  the equations $\tg ij - \gg ij =0$,   define the
Pfaffian system
$$I=\{ \to i - \oo i, \oo\mu, \tooo ij -\ooo ij,
\ooo\mu i-\bb\mu i\oo i \}.\tag\labeleq{\Idef}$$
(There is no sum on $i$ in the last group.)

We will need to compute the derivatives of these forms modulo $I$.  For all
but the last group of forms in (\Idef), the exterior derivatives are zero
modulo $I$.
The exterior derivatives of the forms in the last
group provide the coefficients of  the covariant
derivative of the second fundamental form, $\nabla II$.
Following [BBG], to facilitate computations we write
$$
d\bb\mu i -   \ooo\mu\nu\bb\nu i
 = \sum_j (\bb\mu i -\bb\mu j )\pp ij.
\tag\labeleq{\pidef}
$$
(\pidef) defines the one forms $\pp ij$.
(\pidef) is possible
because it follows from (\bcond) that for any $i$ the vectors
$b_i - b_j$ are a basis for $W$.  By convention, $\pp ii =0$ for any $i$.

It is also convenient to keep track of the difference between the
connections defined by $g$ and $Q$. To this end, write
$$
\tooo ij = \tmu ij + \tnu ij,\tag\labeleq{\wsplit}
$$
where $\tmu ij =\tmu ji$ and $\tnu ij =-\tnu ji$. The $\tmu ij$ measure the
difference between the Levi-Civita connections of $\tilde g$ and $Q$.
Accordingly, there are tensor components
$\tilde \mu^i_{jk}$ such that $\tmu ij=\Sum{k}
\tilde \mu^i_{jk}\to k$ on $\tcf$.  The contraction
$$\nabla Q = (Q_{il}\tilde \mu^l_{jk}+Q_{jl}\tilde \mu^l_{ik})\to i \to j
\otimes \to k$$
is the covariant derivative of $Q$ with respect to the connection of $\tilde
g$.
$\nabla Q$,  denoted by $III$ in [BBG], will play a role in what follows.

There is a (relatively harmless) error in [BBG] in computing
$d(\ooo\mu i-\bb\mu i\oo i )$ modulo $I$.  (See the equation above (4.33) in
that paper.)  We are grateful to Robert Bryant for indicating how to make
a correction.  For the sake of completeness, we include it here:
$$\align
d(\ooo\mu i-
\bb\mu i\oo i )
&
=-\sum_j\ooo\mu j\ww\ooo ji -\sum_\nu\ooo\mu\nu\ww\ooo\nu i
-d\bb\mu i\ww \oo i + \sum_j\bb\mu i\ooo ij\ww\oo j\\
&\equiv  -(d\bb\mu i + \ooo\mu\nu\bb\nu i)
\ww\oo i + \sum_j(\bb\mu i\tooo ij\ww\oo j +
\bb\mu j\tooo ji\ww\oo j) \tmod I
\tag\labeleq{\nicetwo}\\
&\equiv
\sum_j(\bb\mu i -\bb \mu j )\left( -\pp ij\ww\oo i + \tooo ij\ww\oo j \right)
+ 2\sum_j \bb \mu j\tmu ij\ww\to j  \tmod I
\tag\labeleq{\correctvec}
\endalign
$$
(Note that there is no sum on the index $i$.)
Define the functions
$$B_i = b_i \cdot b_i +1$$
on $\Cal W$; then
$$b_i \cdot b_j = -1 + \delta_{ij} B_i\tag\labeleq{\bdot}$$
for any $i,j$.  Dotting (\correctvec) with $b_k$, $k\ne i$, shows that the
  two-forms in the system are
$$B_k( \pp ik \ww \to i -\tnu ik \ww \to k + \tmu ik \ww \to k)
- 2\sum_j \tmu ij \ww \to j,\qquad i \ne k
\tag\labeleq{\corrected}
$$

Compared with [BBG], our formula has an extra term (the one with the
sum over $j$) at the end.  Nevertheless
we will derive the same integrability condition:
$$\tmu ij \equiv 0 \tmod \to i, \to j\ \text{for}\ i\ne j.
\tag\labeleq{\intcondzero}$$
This condition places additional restrictions on the metric $\tilde g$.  In
fact, there
must be contravariant tensor components $\lambda_k$ on $\tcf$ such that
$$\tmu ij = \lambda_i \to j + \lambda_j \to i
   + \delta^i_j \sum_k \lambda_k \to k.\tag\labeleq{\muform}$$
(See pp. 866-867 in [BBG] for details.)  This in turn implies that the extra
term in
(\corrected) is zero, and the rest of the argument in [BBG], from (4.35)
onwards,
goes through with $\tmu ij$ replaced by $-\tmu ij$.

To derive (\intcondzero) from our correction, first note that (\bcond) implies
that
$$\sum_i \dfrac{1}{B_i} = 1.\tag\labeleq{\Bsum}$$
For any $i$ and $j$, let $$\Theta^i_j = \too ij \ww \to j \ww \to i.$$
Then wedging the two-form
(\nicetwo) with $\to i$ gives $\Sum{j} (b_i \Theta^i_j - b_j \Theta^j_i)$.
Dotting with $b_k$ for $k\ne i$ gives the three-form $-B_k\Theta^k_i
+ \Sum{j}(\Theta^j_i - \Theta^i_j)$.  Thus, on any integral manifold of $I$,
$$\Theta^k_i = \Sum{j}(\Theta^j_i - \Theta^i_j)/B_k.\tag\labeleq{\krusher}$$
Summing over $k\ne i$, using (\Bsum) and solving gives
$$\sum_j\Theta^i_j = (1+B_i)\sum_j\Theta^j_i.$$
Substituting this into (\krusher) and summing both sides over $k$ gives
$$\sum_j\Theta^j_i = (1-B_i)\sum_j\Theta^j_i.$$
Since $B_i \ne 0$, it follows that, on any integral manifold, $\Theta^i_j=0$
for every $i$ and $j$.  The integrability condition now follows
using the decomposition (\wsplit).

Having shown that there exist quasi-hyperbolic metrics for which this
integrability condition holds, [BBG] go on to determine whether or not the
system is involutive.
Its tableau has characters
$s_1=n(n-1)$, $s_2=\hdots = s_n=0$.  Integral elements are
obtained by setting all forms in $I$ equal to zero, choosing $n(n-1)$ constants
$\{ A^i_j | i\neq j \}$ and setting
$$
\aligned
\tnu ij &=-\lambda_i\to j + \lambda_j\to i
+A^i_jB_j\to j -A^j_iB_i\to i \\
\pi_{ij} &= A^j_iB_i\to j - A^i_jB_i\to i.
\endaligned \tag\labeleq{\element}
$$
Since the space of integral elements is $n(n-1)$ dimensional,
the system is involutive and local solutions depend
on $n(n-1)$ functions of one variable.

\bigpagebreak
\heading \S3. The minimality condition \endheading
Now  we add the requirement that the image of the
isometric embedding be a minimal submanifold.   The minimality condition
$\Sum{i,j}g^{ij} h^\mu_{ij}=0$ and (\hform) imply that
$$
\sum_i\bb{} i\ugg ii =0. \tag\labeleq{\firstmincond}
$$
Our system for minimal isometric immersions
will be $I$ restricted to the submanifold $\Sigma'''
\subset
\Sigma''$ where (\firstmincond) holds.

Any set of vectors $b_i$ satisfying (\bcond) have exactly one linear relation
among them, and up to multiple this must be
$$
\sum_i \dfrac{b_i}{B_i} = 0.\tag\labeleq{\tidbit}
$$
Using (\tidbit)
 and (\Bsum), one gets an equivalent minimality
condition,
$$\dfrac{1}{B_i} = \dfrac{\ugg ii}{\Sum{j} \ugg jj}.
\tag\labeleq{\secondmincond}$$
In order to see if any integral elements now exist for $I$, we will need
to see if (\element) is compatible with the additional linear relations on
the one-forms $\pi_{ij}$ and $\too ij$ introduced by the restriction to
$\Sigma'''$.

On $\Sigma''$, the relations (\bcond) implied   $dB_i = 2B_i
\Sum{j}\pi_{ij}$.
We already had the relations
$$
B_j\pi_{ij} + B_i\pi_{ji}=0,\qquad i\ne j.\tag\labeleq{\pirel}
$$
We also note that
$$d\ugg ij = -\ugg ik \too jk - \ugg jk \too ik.$$
Now, differentiating (\secondmincond) gives
$$\sum_j \dfrac{\pi_{ij}}{B_i}
  = (\sum_j \ugg ij \too ij - \ugg ii\sum_{j,k}\ugg jk\too jk)/(\Sum{j}
\ugg jj).
\tag\labeleq{\stepone}$$
These constitute $n-1$ additional linearly independent relations on the
$\pi_{ij}$.
Substituting (\muform) and (\element) into (\stepone) gives
$$B_i\ugg ii \sum_j(A^j_i \to j-A^i_j \to i)=
     \sum_j \ugg ij(A^i_j B_j \to j - A^j_i B_i\to i+2\lambda_i\to j)
    -2\dfrac{1}{B_i}\sum_{j,k}\ugg jk\lambda_j \to k,$$
where for convenience we set $A^i_i = 0$ for all $i$.  For $k\neq i$, taking
the coefficient of $\to k$ on each side gives
$$B_i \ugg ii A^k_i =\ugg ik (A^i_k B_k + 2 \lambda_i)
            -2\dfrac{1}{B_i}\sum_j \ugg jk\lambda_j.$$
(The equations obtained by taking the coefficient of $\to i$ will be linearly
dependent on the last set.)
This equation, together with that obtained by interchanging $i$ and $k$,
gives a pair of linear equations for $A^k_i$ and $A^i_k$,
$$
\aligned B_i \ugg ii A^k_i - B_k \ugg ik A^i_k
&= 2 \ugg ik\lambda_i -2\dfrac{1}{B_i}\sum_j\ugg jk\lambda_j\\
 B_k \ugg kk A^i_k -B_i \ugg ik A^k_i
&= 2 g^{ik}\lambda_k -2\dfrac{1}{B_k}\sum_j \ugg ij\lambda_j.\endaligned
\tag\labeleq{\Aeqns}
$$
Since the coefficient matrix on the left has determinant
$B_i B_k (\ugg ii \ugg kk - (\ugg ik)^2)$, and this is clearly nonzero, these
 equations determine the
coefficients
$A^i_j$  uniquely in terms of functions defined on $\tcf$.  We conclude that
there is a unique integral $n$-plane satisfying the independence condition at
each point of $\Sigma'''$.  In fact, we may define a new Pfaffian system
on $\Sigma'''$ by adding more generator 1-forms to $I$:
$$J =  I \oplus \{ \tnu ij +\lambda_i\to j - \lambda_j\to i
-A^i_jB_j\to j +A^j_iB_i\to i,
\pi_{ij} - A^j_iB_i\to j + A^i_jB_i\to i\},\tag\labeleq{\defJ}$$
where the $A^i_j$ are determined by (\Aeqns).
At each point of $\Sigma '''$, $J$ annihilates our
distribution of $n$-planes.  Any integral manifold of $I$ in $\Sigma'''$
satisfying
the independence condition will be an integral manifold of $J$.
At this point we see that solutions depend at most on constants.

\demo{Proof of Theorem C}  When $(\tm,\tilde g)$ is itself a space form
  of constant sectional curvature $k$, a trace of the Gauss equations implies
(as in the proof of [M], Thm. 2) that $\epsilon \ge k$.  If $\epsilon =k$ then
$M$ is totally geodesic.  If $\epsilon > k$, the Gauss equations take the
form (\flatgauss) with $\widehat{II} =II \oplus \sqrt{\epsilon - k}\, \tilde
g$.
Then Cartan's theorem implies the rank is at least $n-1$.

Now assume the rank is $n-1$; if the codimension exceeds the rank, extra
1-forms
are present in the system $I$, but the 2-forms are unchanged, so we add the
same 1-forms as in (\defJ) to obtain system $J$. When we take
$Q=\sqrt{\epsilon -k}\,\tilde g$, the forms $\tmu ij$ are automatically zero,
as are the right-hand sides in (\Aeqns).  Now $J$ takes the form
$$J =  I \oplus   \{\tnu ij,\pi_{ij}\}.$$
However, because
$$\align d\tnu ij &= -\tnu ik \ww \tnu kj - \to i \ww \to j \\
&\equiv - k\to i \ww \to j\ \tmod J,\endalign$$
we see that unless $k=0$,
$J$ satisfies the Frobenius condition nowhere on $\Sigma'''$.
Uniqueness in the flat case now follows by the argument
given at the end of the proof of Theorem D.

More generally, if $Q$ is $\tilde g$-parallel, then the forms $\tmu ij$
and tensor components $\lambda_i$  must vanish.  Once again,
$J =  I \oplus   \{\tnu ij,\pi_{ij}\}$, and the argument proceeds as above.
\enddemo

\smallpagebreak

Returning to the general case, we assume now that $(\tm^n,
\tilde g)$ is a quasi-$\k$ curved Riemannian manifold with respect to a
positive-definite $Q$, such that (\muform) holds.

Any integral $n$-manifold of $J$ will push down to an integral $n$-manifold
in $\tcf$ for the system
$$\tilde J = \{\tnu ij +\lambda_i\to j - \lambda_j\to i
-A^i_jB_j\to j +A^j_iB_i\to i\},$$
where $B_i$ and $A^i_j$ now are determined by
by the functions
$\tilde \ugg ij$ and $\lambda_i$ on
$\tcf$.
So, the half of the Frobenius condition for $J$ that involves
$d\tnu ij$ will only require the vanishing of certain functions on $\tcf_Q$,
and these are intrinsic integrability conditions that
depend only on the metric $\tilde g$ and on $Q$.
In order to examine the other half, we need to compute $d\pp ij$ modulo $I$.

Differentiating (\pidef) gives
$$
\multline
\sum_{j,k} \bb \mu j \tilde g^{jk} (b_i \cdot b_k) \to j \ww \to k \\
\equiv \sum_j \left[ \left( \sum_k (\bb \mu i - \bb \mu k) \pp ik -
\sum_k (\bb \mu j - \bb \mu k) \pp jk\right) \ww \pp ij + (\bb \mu i - \bb
\mu j)
d\pp ij\right]
\tmod I.\endmultline
$$
Dotting with $b_p$, $p\ne i$, gives
$$
\multline
\sum_{j,k}\tilde g^{jk} (b_j \cdot b_p) (b_i \cdot b_k) \to j \ww \to k \\
\equiv \sum_j \left( B_p \pp ip -
\sum_k ((b_j \cdot b_p) - (b_k \cdot b_p)) \pp jk\right) \ww \pp ij +
B_p d\pp ip
\tmod I.\endmultline
$$
Using (\bdot), we see that the
value of $d\pp ij$ modulo $I$ can be expressed in terms of the $\pp ij$'s
themselves and functions and forms defined on $\tcf_Q$,
so these integrability conditions also are intrinsic.
In other words,
{\it the Frobenius condition for $J$ can be expressed solely in terms of the
vanishing of certain functions on $\tcf_Q$.}

\proclaim{Proposition} Given a section $\{\tilde e_i\}$ of $\tcf_Q$ along
which the aforementioned functions vanish, there exists a minimal isometric
embedding of $M$ into $X$, unique up to rigid motion, under which $\{\tilde
e_i\}$
are the principal tangent directions.\endproclaim
\demo{Proof} Existence follows from the Frobenius theorem.  Suppose now there
are two such embeddings $f,F$.  Fixing a point
$p\in M$, we can arrange
by rigid motions that $f(p)=F(p)=x$ and $f_*(T_pM) = F_*(T_pM)=T_x$.  Now
we have
$$g_{ij}= <f_*\tee i,f_*\tee j>=<F_*\tee i,F_*\tee j>$$
and it follows that we can arrange, by rotations of $X$ acting on the
plane $T_x$, that $f_*\tee i=F_*\tee j$.  We can also arrange that
the principal normal vectors
$$b_i = \sum_\mu b^\mu_i e_\mu\in N_x$$ also coincide.
(For, the $b_i$ are $n$ vectors with
$b_i \cdot b_j = -1 +\delta_{ij}\sum_k \tilde g^{kk} /\tilde g^{ii}$, and
the set of such vectors is acted on simply transitively by rotations in $N_x$.)
To each of $f$ and $F$ there is associated
a ``graph'', a section of $\Sigma''' \subset \tcf_Q\times\cf\times \Cal W$
which covers the embedding, and which is an integral manifold of $J$.  By using
the action of $O(n-1)$ along the Cauchy characteristics of $J$, i.e. the action
by
orthogonal substitions among the frame vectors $e_\mu$, we can arrange that the
frames associated to $f$ and $F$ coincide at $x$.  Then the vectors $b^\mu_i$
must coincide there as well.  This means that the integral manifolds of $J$
associated to $f$ and $F$ go through the same point of $\Sigma'''$ above $p$,
and so the embeddings must coincide everywhere else.
\enddemo

\demo{Proof of Theorem D}
Again, fix a generic point $p \in M$, and let $G$ be the identity component
of the group of isometries of $M$ fixing $p$. Then, because $G$ fixes the Ricci
tensor at $p$, it also fixes $Q$.
The argument of the above proposition could be used to prove that a minimal
isometric embedding of $M$ is unique up to rigid motion if we knew that
$G$ acted simply transitively on $Q$-orthonormal frames $\{\tee i\}$.
(This is the case, of course, when $Q=\tilde g$.)   For, if $f$ and $F$ are two
such embeddings, and the associated sections of $\tcf_Q$ are, at $p$,
in the same orbit of $G$, then we can arrange by isometries of $M$ that
they are the same at $p$, and (as in the proposition) arrange by rigid motions
of $X$ that the embeddings coincide.  Hence the embeddings depend, up to
rigid motions, only on the value of $\{\tee i\}$ at $p$, modulo the action of
$G$. \enddemo

\heading \S4. Isometric Embeddings with Parallel Mean Curvature \endheading
Suppose $(\tm,\tilde g)$ is a quasi-$\k$-curved Riemannian manifold, and $\tm$
is isometrically embedded in space form $X^{2n-1}(\k+1)$.
In terms of the usual adapted framing (see \S 2), the mean curvature vector is
$$H = \sum_{\mu,i,j} e_\mu g^{ij} h^\mu_{ij} = \sum_{\mu,i} e_\mu g^{ii}b^\mu_i
.$$
Using (\pidef) and the structure equations of $\cf$, the normal component of
$\nabla H$ is
$$\nabla_N H = \sum_{\mu,i,j} e_\mu\left(g^{ii}(b^\mu_i-b^\mu_j)\pi_{ij}
-2g^{ij}b^\mu_i  \omega^i_j\right).$$
Thus, requiring that the mean curvature vector be parallel amounts to requiring
that the vector-valued 1-forms
$$\sum_{i,j} \left(g^{ii}(b_i-b_j)\pi_{ij}
-2g^{ij}b_i  \omega^i_j\right)$$
vanish, in addition to the forms in the system $I$ defined by (\Idef).  Dotting
with
vector $b_k$ gives the equivalent condition
$$0=\sum_j \left(\pi_{kj}g^{kk}-\pi_{jk} g^{jj} - 2g^{jk} \too kj\right)
 + {2 \over B_k}\sum_{i,j}g^{ij}\too ij.\tag\labeleq{\mcforms}$$
Thus, such isometric embeddings will arise as integral manifolds of the
following Pfaffian
system on $\Sigma''$:
$$K = I \oplus \{\sum_j \left(\pi_{kj}g^{kk}-\pi_{jk} g^{jj} - 2g^{jk} \too
kj\right)
 + {2 \over B_k}\sum_{i,j}g^{ij}\too ij\}.$$

Now we may substitute in (\mcforms) the values of $\pi_{ij}$ and $\too ij=\tmu
ij + \tnu ij$ on a typical integral
element of $I$, as given by (\muform) and (\element), to get the condition
$$\multline
\sum_jA^j_k(B_k g^{kk}+B_j g^{jj})\to j - A^k_j(B_k g^{kk}+B_j g^{jj})\to k
-2g^{jk}(A^k_jB_j \to j -A^j_k B_k \to k) \\
=4(\Sum{j} g^{jk}\lambda_j \to k- {1\over B_k}\Sum{i,j} g^{ij}\lambda_j \to i)
+2(g^{kk} - {1\over B_k} \Sum j g^{jj}) (\Sum l\lambda_l \to l)
\endmultline$$
For $i\ne k$, equating the coefficients of $\to i$ on both sides gives
$$\align
(B_k g^{kk}+B_i g^{ii})A^i_k -2g^{ik}B_i A^k_i
&= 2\lambda_i(g^{kk} - {1\over B_k} \Sum jg^{jj})
-{4\over B_k} \Sum j g^{ij}\lambda_j\\
-2g^{ik}B_k A^i_k + (B_k g^{kk}+B_i g^{ii})A^k_i
&= 2\lambda_k( g^{ii} - {1\over B_i} \Sum j g^{jj})
-{4\over B_i} \Sum j g^{kj}\lambda_j
\endalign$$
(The second equation comes from the first by interchanging $i$ and $k$.)  This
gives two linear equations for $A^i_k$ and $A^k_i$, and the determinant
of the coefficient matrix is
$$(B_k g^{kk}+B_i g^{ii})^2 - 4B_i B_k (g^{ik})^2
=(B_k g^{kk}-B_i g^{ii})^2 + 4B_i B_k (g^{ii}g^{kk} -(g^{ik})^2) > 0.$$
So, as happened with the minimality condition in \S 3, at each point of
$\Sigma''$
there is a unique integral $n$-plane for $K$.  Now the arguments in the proofs
of Theorems C and D apply as before.  (Note, however, that tracing the
Gauss equations, as we do in the proof of Theorem C, does not yield any useful
inequalities for the curvature of $M$ and $X$ in this case.)

\heading \S 5. Construction of quasi-$\k$-curved metrics \endheading

Let $\ee 0\hd \ee{n+1}$ be an orthornormal frame of $\Bbb L^{n+2}$
with $\ee 0$ a timelike direction.
On the orthonormal frame bundle $\cf_{\Bbb L}$
of $\Bbb L^{n+2}$,
 we have the Maurer-Cartan form
$$
\Omega =
\pmatrix
0&0&0&0\\
\oo 0 &0&\ooo 0i & \ooo 0{n+1}\\
\oo j & \ooo j0 & \ooo ij & \ooo j{n+1}\\
\oo{n+1} & \ooo{n+1} 0&\ooo{n+1}j & 0
\endpmatrix
$$
where $\ooo ij=-\ooo ji, \ooo 0{i}=\ooo{i}0,
\ooo 0{n+1}=\ooo{n+1}0,\ooo j{n+1}=-\ooo{n+1}j$.

Assume $M^n\subset\bll{n+2}$ is spacelike.
Let $\cf^1_{\bll{}}$
denote the subbundle of
 the restriction of
$\cf_{\bll{}}$ to $M$  consisting of frames such that
 at each point $x\in M$, $T_xM=\{\ee 1\hd\ee{n}\}$. On
$\cf^1_{\bll{}}$, one has
$$
\ooo 0 i = q_{ij}\oo j, \ \ \ooo{n+1}i=h_{ij}\oo j
$$
for some functions $q_{ij}, h_{ij}$ symmetric in their
lower indices. The curvature tensor of $M$ is
$$
R= \gamma (h,h)-\gamma (q,q)
$$
To construct quasi-$\k$-curved metrics, we want
$$
h=\sqrt{\k+1}g
$$
where $g$ is the induced metric on $M$.
Since we are working with $g$-orthonormal frames,
we need $h_{ij}= \sqrt{\k+1}\delta_{ij}$.

In what follows, we will repeatedly differentiate the above conditions,
 using the Maurer-Cartan structure equations.  Since we will be using the
results
of [Lan], we will generally use formulae and notation from there,
except that the indices $i,j,k$ will index the tangent directions, instead
of $\alpha, \beta, \gamma$. $\mu,\nu$ will index normal directions, in this
case just $0$ and $n+1$.  In particular, $\rr\mu ijk$ will
denote the coefficients of $F_3$, which in our situation
 is a tensor on $M$, namely $\nabla II$.
Similarly, $F_k = \nabla^{k-2} II$.

The coefficients of $F_3,F_4,F_5$ are obtained by differentiating
$\ooo\mu j-\qq\mu jk\oo k$. They are given as follows:
$$
\align
\rr\mu ijk\oo k
&=
-d\qq\mu ij -\qq\nu ij\ooo\mu\nu + \qq\mu il\ooo l j + \qq\mu jl\ooo li
\\
\rr\mu ij{kl}\oo l
&=
-d\rr\mu ijk -\rr\nu ijk\ooo\mu\nu +
\frak S_{ijk}\rr\mu ijl\ooo lk -
\frak S_{ijk}\qq\mu il\qq\nu jk\ooo l\nu
\\
\rr\mu ij{klm}\oo m
&=
-d\rr\mu ij{kl} - \rr\nu ij{kl}\ooo\mu\nu
+\frak S_{ijkl}\rr\mu ij{km}\ooo ml
-\frak S_{ijkl}(\rr\mu ijm\qq\nu kl + \qq\mu im \rr\mu jkl )\ooo m\nu,
\endalign
$$
where $\frak S_{ijk}$ denotes a cyclic sum over $i,j,k$ and
$\frak S\frak S_{ijkl}$ denotes summing to symmetrize the expression
over $i,j,k,l$.

Finally, because results in [Lan] are phrased in terms of
submanifolds of complex projective space, we will
consider $\bll{n+2}\subset\brrr\bpp{n+2}\subset\bccc\bpp{n+2}$ and
$\cf_{\bll{}}\subset\cf_{\bccc\bpp{n+2}}$.

\demo{Proof of Lemma}
To indicate how the proof applies to more general signatures, we will let
$\epsilon_{\mu}= - <\ee\mu ,\ee\mu >=\pm 1$;
then $\ooo \mu j=\epsilon_\mu \ooo j\mu$.

Using ([Lan], 2.15), or by differentiating the equation $\ooo{n+1}j-\sqrt{\k+1}
\oo j=0$,
we obtain
$$
\rr{n+1}ijk\oo k = -q_{ij}\ooo{n+1}0
\qquad\text{and}\qquad
\ooo{n+1}0\wedge \epsilon_0 q_{ij}\oo j=0\ \ \forall i.
$$
Since $Q$ is nondegenerate, this implies that
$$
\rr{n+1}ijk =0\text{ and }\ooo{n+1}0 =0
$$
which proves the assertion $\sigma\intprod \nabla II=0$.
(Note that  the normal bundle splits into
parallel sub-bundles spanned by $\sigma$ and $\sigma^{\perp}$.)  Now
$$
\align
\rr{n+1}ij{kl}\oo l
&= \frak S_{ijk}\qq{n+1}im\qq\mu jk\ooo m\mu \\
&= \epsilon_{\mu}\sqrt{\k+1}\frak S_{ijk}\qq{\mu}im\qq\mu jk\oo m  \endalign
$$
so that
$$
\rr{n+1}ij{kl}=
\epsilon_{\mu}\sqrt{\k+1}\frak S\frak S_{ijkl}\qq{\mu}il\qq\mu jk,
$$
which implies that ([Lan], 4.13) holds with
$b^{n+1}_{\mu\nu}=\delta_{\mu\nu}\epsilon_{\mu}\sqrt{\k+1}$.

Now to compute coefficients of $F_5$:
$$
\align
\rr{n+1}ij{klm}\oo m
&= -d\rr{n+1}ij{kl} + \frak S_{ijkl}\rr{n+1}ij{km}\ooo ml
-\frak S_{ijkl}\qq{n+1}im\rr\mu jkl\ooo m\mu \\
&=
\epsilon_{\mu}\sqrt{\k+1} \frak S\frak S_{ijkl}
d(\qq\mu ij\qq\mu kl ) -
\frak S_{ijk}(\frak S\frak S_{ijkm}\qq\mu ijk\qq\mu km)\ooo ml )
-\epsilon_{\mu}\frak S_{ijkl} \rr\mu jkl\qq\mu im\oo m
\endalign
$$
When we substitute
$$
d\qq\mu ij = - \rr\mu ijm\oo m + \qq\mu im\ooo mj + \qq\mu jm\ooo mi
-\qq\nu ij\ooo\mu\nu,
$$
everything in $\rr{n+1}ij{klm}\oo m$ cancels except a term of the form
$\epsilon_{\mu}\rr\mu ijk\qq\mu lm\oo m$, which implies that
([Lan], 4.16) holds as well.
This implies that at each point there is a quadric, whose tangent space
is $\ee{n+1}\upperp$, osculating to order five.
Since we are in codimension two, $Q$ is nonzero, and there is at least one
quadric of   rank at least three,
we know there are no linear syzygies among the quadrics in $II$. We thus
  conclude by ([Lan], 4.20) that the quadric actually contains $M$.
Computing at the point $(0\hd 0)$, following ([Lan], 4.18)
we find that the equation of the quadric is:
$$
\xx{n+1}  - \sqrt{k+1}((x^1)^2 + \cdots +(\xx n)^2 -(\xx 0)^2)
+(x^{n+1} )^2)
=0
$$
and completing the square yields (\quadric).
\enddemo

We now examine the quasi-$\k$-curved metrics satisfying the
integrability conditions; i.e., we want metrics such that
$$
\ee 0\intprod \nabla II = L\circ (\ee 0\intprod II). \tag\labeleq{\nabii}
$$
This implies that at each point of $M$ there
is a quadric hypersurface
having tangent plane $\ee 0\upperp$   osculating to order three.
  Having a quadric osculate to order
three at a general point is not enough to imply containment;
but our additional conditions
will imply containment.

We will generally follow [Lan]  in the remainder of this section
with one exception:
in what follows $0$ will denote a normal index, and what is denoted
$\xx 0$ in [Lan]  should be taken to be a homogeneous coordinate which
 we set equal to $1$ in restricting to an affine space.

What follows is actually a result in projective geometry; only one
needs the hypotheses that one has two sections $\sigma$ and $\tau$
of the normal bundle $NM$ which have the property that
$\sigma\intprod F_3=0$ and
$\tau\intprod F_3= L\circ (\tau\intprod II)$. Then one can restrict to
$\sigma\intprod II$-orthonormal frames in $TM$ and
to frames in $NM$ where $\sigma$ and $\tau$ form
pointwise a Minkowski orthonormal basis, and then restrict
further to an affine open subset.
For simplicity, we will work in Lorentz space.

In frames, our hypothesis is that
$$
\rr 0ijk =\frak S_{ijk}\lm i\qq 0 jk \tag\labeleq{\rcoefs}
$$
i.e. that
$$
-d\qq 0 ij + \qq 0 il\ooo lj + \qq 0 jl\ooo li
= \Sum k (\lm i\qq 0 jk+\lm j\qq 0 ik+\lm k\qq 0 ij)\oo k
 \tag\labeleq{\fthreeb}
$$

Computing the coefficients of $F_4$ (i.e. differentiating
(\rcoefs)) we obtain
$$
\align
\rr 0 ij{kl}\oo l
&=
-d\rr 0 ijk
+\frak S_{ijk} \rr 0 ijm\ooo mk
-\frak S_{ijk}\qq 0 im\qq 0 jk\ooo m0
 -\frak S_{ijk}\qq 0 im\qq {n+1} jk\ooo m{n+1}\\
&=
\frak S_{ijk} \left\{
-\qq 0 jkd\lambda_i + \lambda_l\qq 0 ij\ooo lk +
 [
\lambda_i(\frak S_{jkl}
\lambda_j\qq 0 kl) + (\k+1)\qq 0 il\delta_{jk}-\Sum m\qq 0 im\qq 0 ml\qq 0 jk
 ]\oo l \right\}
 \tag\labeleq{\ffoura}\endalign
$$

  Note that (\ffoura)
includes a term $LF_3$, which is necessary for the
Monge system, a term of   potential torsion (obstruction to
integrability), and a term of
Cauchy characteristics.

To simplify computations, at this point we
restrict to frames where $q^0$ is diagonal and write
$q_i$ for $q^0_{ii}$. Then (\fthreeb) becomes
$$
\align
&(q_i-q_j)\ooo ij =
\lambda_iq_j\oo j + \lambda_j q_i\oo i \tag\labeleq{\fthreec}\\
& d(log (q_i))=
-\Sum k\lambda_k\oo k - 2\lambda_i\oo i. \tag\labeleq{\fthreed}\endalign
$$
If $q_i\neq q_j$ we get
$$
\ooo ij =\frac{\lambda_iq_j}{q_i-q_j}\oo j
+\frac{\lambda_jq_i}{q_i-q_j}\oo i \tag\labeleq{\fthreee}
$$
If none of $i,j,k$ are equal, (\ffoura) simplifies to
$$
\rr 0 ij{kl}\oo l =
2q_k\lambda_i\lambda_j\oo k +
2q_j\lambda_k\lambda_i\oo j +
2q_i\lambda_j\lambda_k\oo i \tag\labeleq{\ffourb}
$$
Note that this is the fourth order Monge condition ([Lan], 4.17)
for these indices.

For simplicity, we first assume that the eigenvalues $q_i$ of $Q$ are distinct
on an open set in $M$.
(\ffourb) implies that
$$
\align
&\rr 0 ij{kl}=0\hphantom{\lambda_j\lambda_kq_i }
 \text{if }ijkl\text{ are all distinct,}\\
&\rr 0 ii{jk}= 2\lambda_j\lambda_kq_i
\hphantom{0 }
\text{if }ijk\text{ are all distinct.}
\tag\labeleq{\rzeros}
\endalign
$$
On the other hand, when $i\ne k$, from (\ffoura) we get
$$
\align
\rr 0 ii {kl}\oo l = -&q_id\lambda_k +
(4\lambda_i\lambda_kq_i)\oo i \\
&+[2\lambda_i^2q_k + \lambda_k^2q_i + (\k+1)q_k - q_k^2q_i +
q_iq_k(\Sum{l\neq k}\frac{\lambda_l^2}{q_l-q_k})]\oo k\\
&+\Sum{l\neq k}[\lambda_l\lambda_kq_i(\frac{q_l}{q_l-q_k}+1)]\oo l
\tag\labeleq{\fthreef}
\endalign
$$
The symmetry of $\rr 0 ij{kl}$ in the
lower indicies places restrictions on $d\lambda_i$. If we let
$d\lambda_k= \lambda_{kl}\oo l$,
then (\fthreef) implies, for $i,k,l$ distinct,
$$
\align &
\rr 0 ii{kl} = -q_i\lambda_{kl} + \lambda_k\lambda_lq_i(
\frac{q_l}{q_l-q_k} + 1)\tag\labeleq{\fthreeg}\\
&\rr 0ii{kk} = q_i\left(-\lambda_{kk}+\lambda_k^2-q_k^2
+q_k(\Sum{l\neq k}\frac{\lambda_l^2}{q_l-q_k})\right)
 +(\k+1)q_k+ 2\lambda_i^2q_k
\tag\labeleq{\fthreeh}
\endalign
$$

Combining (\rzeros) and (\fthreeg) we obtain
$$
\lambda_{kl}= \frac{\lambda_k\lambda_lq_k}{q_l-q_k}, \ \ k\neq l.
$$
Using the symmetry $\rr 0 ii{kk}=\rr 0{kk}ii$, we obtain
$$
q_i(\lambda_{kk} +\lambda_k^2+(\k+1)+q_k^2
-q_k\Sum{l\neq k}\frac{\lambda_l^2}{q_l-q_k})=
q_k(\lambda_{ii} +\lambda_i^2+(\k+1)+q_i^2
-q_i\Sum{l\neq i}\frac{\lambda_l^2}{q_l-q_i}).
$$
Hence
$$
\frac 1{q_k}(\lambda_{kk} +\lambda_k^2+(\k+1)+q_k^2
-q_k\Sum{l\neq k}\frac{\lambda_l^2}{q_l-q_k}) \tag\labeleq{\fthreei}
$$
is independent of $k$; call this quantity $-b^0_{00}$
and set $b^0_{0,n+1}=\k+1$.  Then
$$
\rr 0 ii{kk}= 2\lambda_i^2q_k + 2\lambda_k^2q_i +b^0_{00}q_iq_k +
b^0_{0,n+1}(q_i+q_k),
$$
and thus the remaining fourth order Monge conditions
in the $e_0$ direction hold as well:
$$
F^0_4= LF^0_3 + b^0_{00}F^0_2F^0_2+b^0_{0,n+1}F^0_2F^{n+1}_2.
\tag\labeleq{\mongefour}$$
Differentiating again, one sees that
the fifth order Monge condition
 ([Lan],4.18)
holds as well.
Setting $b=b^0_{00}$ and using ([Lan], 4.18) one gets (\tquadric).

\smallpagebreak

Now assume that some eigenvalues of $Q$ have multiplicity one and some have
multiplicity greater than one on an open set in $M$; let $\xi,\eta$ index
the former and $\alpha,\beta$ index the latter.  Then (\fthreec) implies
$\lambda_\alpha=0$ and
$$\ooo \xi\alpha =\frac{\lambda_\xi q_\alpha}{q_\xi-q_\alpha}\oo \alpha.$$
Differentiating this gives the analogous results:
$$
\lambda_{\xi\alpha}=0,\qquad
\lambda_{\xi\eta}= \frac{\lambda_\xi\lambda_\eta q_\xi}{q_\eta-q_\xi},
\ \ \xi\neq \eta,
$$
and
$$
\frac 1{q_\xi}(\lambda_{\xi\xi} +\lambda_\xi^2+(\k+1)+q_\xi^2
-q_\xi\Sum{\eta\neq \xi}\frac{\lambda_\eta^2}{q_\eta-q_\xi})
=
\frac 1{q_\alpha}((\k+1)+q_\alpha^2
-q_\alpha\Sum{\eta}\frac{\lambda_\eta^2}{q_\eta-q_\alpha})
= -b^0_{00},
$$
independent of $\alpha$ and $\xi$.  Now the fourth-order Monge condition
(\mongefour) holds without change, as does the fifth-order condition.

\smallpagebreak

The case where the eigenvalues of $Q$ all have multiplicity greater than one is
simple: (\fthreec) implies that $\lambda_i=0$ for all $i$, hence the $q_i$ are
constant, and $M$ is a product of space forms given
by the intersection of (\quadric) and  (\zquadric).

\heading Acknowledgements \endheading
We would like to thank Robert L. Bryant for helpful discussions and advice,
and to thank the referee for several comments and suggestions.

\heading References \endheading
\parindent=0pt
[BBG] E. Berger, R. Bryant, P. Griffiths, {\it The Gauss equations and
rigidity of
isometric embeddings}, Duke Math. J. {\bf 50} (1983) 803-892.
\smallskip
[Cal] E. Calabi, {\it Isometric imbeddings of complex manifolds},
 Annals of Math.  {\bf 58} (1953) 1-23.
\smallskip
[Car1] E. Cartan, {\it Sur les vari\'et\'es de courbure constante d'un
espace euclidien ou non euclidien},
Bull. Soc. Math France. {\bf 47} (1919) 125-160
and  {\bf 48} (1920), 132-208;
Oeuvres Compl\`etes (Gauthier-Villars, 1955), Ptie. 3, Vol. 1, 321-432.
\smallskip
[Car2] E. Cartan, {\it Les syst\`emes differ\'entielles ext\'erieurs et leurs
applications g\'eo\-m\'et\-riques}, Hermann (1945).
\smallskip
[CO] S.S. Chern, R. Osserman, {\it Remarks on the Riemannian
metric of a minimal submanifold}, Springer LNM. {\bf 894} (1980) 49-90.
\smallskip
 [GH] P.A. Griffiths,   J. Harris, {\it Algebraic Geometry and Local
Differential Geometry}, Ann. scient. Ec. Norm. Sup.  {\bf 12} (1979), pp.
355-432.
\smallskip
[Lan] J.M. Landsberg {\it Differential-geometric characterizations
of complete intersections}, J. Differential Geometry {\bf 44}(1996), pp.
32-73.
\smallskip
[Law] H.B. Lawson, {\it Lectures on
minimal submanifolds, Vol. 1}, p. 178ff; Publish or Perish, 1980.
\smallskip
[M] J.D. Moore, {\it Isometric immersions of space forms in space forms},
Pacific J. Math  {\bf 40}(1972) 157-166.
\smallskip
[T] T.Y. Thomas, {\it Riemann spaces of class one and
their characterization}, Acta Math.   {\bf 67}(1936)
169-211.
\enddocument